\documentclass[submission%
,nohyperref%
]{dmtcs}

\usepackage[latin1]{inputenc}
\usepackage{subfigure}
\usepackage{amsmath}
\usepackage{amssymb}

\usepackage{multirow}

\usepackage{float}

\usepackage{graphicx}
\usepackage{calc}
\newlength{\depthofsumsign}
\newtheorem{defi}{Definition}[section]
\newtheorem{teor}{Theorem}[section]

\newtheorem{lemma}{Lemma}[section]
\newtheorem{conj}{Conjecture}[section]
\newtheorem{corollary}{Corollary}[section]

\newtheorem{remark}{Remark}[section]

\newcommand{\hector}[1][1.65]{
    \mathop{%
        \raisebox
            {-#1\depthofsumsign+1\depthofsumsign}
            {\scalebox
                {#1}
                {$\chi$}%
            }
    }
}

\newcommand{\chisotazo}[1][3pt]{%
  \mathrel{\raisebox{#1}{$\hector$}}%
}

\usepackage[round,numbers]{natbib}

\author{H\'ector Blandin \addressmark{1}\thanks{Email: \email{hectorblandin@gmail.com}. Partially supported by ISM and LaCIM}}
\title{Generalized Polarization Modules (extended abstract)}

\address{Laboratoire de Combinatoire et d'Informatique Math\'ematique (LaCIM), UQ\`{A}M, Canada\addressmark{1}}


\received{20...-??-??}
\revised{\today}
\accepted{?????}

\begin{document}

\maketitle

\begin{abstract}

\paragraph{Abstract.} This work enrols the research line of M. Haiman on the Operator Theorem (the old operator conjecture). This theorem states that the smallest $\mathfrak{S}_n$-module closed under taking partial derivatives and closed under the action of polarization operators that contains the Vandermonde determinant is the space of diagonal harmonics polynomials. We start generalizing the context of this theorem to the context of polynomials in $\ell$ sets of $n$ variables $x_{ij}$ with $1\leq i\leq \ell$ et $1\leq j\leq n$. Given a $\frak{S}_n$-stable family of homogeneous polynomials in the variables $x_{ij}$ the smallest vector space closed under taking partial derivatives and closed under the action of polarization operators that contains $F$ is the polarization module generated by the family $F$. These polarization modules are all representation of the direct product $\mathfrak{S}_n\times{GL}_{\ell}(\mathbb{C})$. In order to study the decomposition into irreducible submodules, we compute the graded Frobenius characteristic of these modules. For several cases of $\mathfrak{S}_n$-stable families of homogeneous polynomials in $n$ variables, for every $n\geq 1$, we show general formulas for this graded characteristic in a global manner, independent of the value of $\ell$.

\paragraph{R\'esum\'e.} Ce travail s'inscrit dans la lign\'ee de recherche des travaux de M. Haiman sur le th\'eor\`{e}me de l'op\'erateur (ex-conjecture de l'op\'erateur). Ce th\'eor\`{e}me affirme que le plus petit $\frak{S}_n$-module clos par d\'erivation partielle et clos par l'action des op\'erateurs de polarisation qui contient le d\'eterminant de Vandermonde est l'espace des polyn\^omes harmoniques diagonaux. On commence par g\'en\'eraliser le contexte du th\'eor\`eme de l'op\'erateur au contexte de polyn\^omes \`{a} $\ell$ ensembles de $n$ variables $x_{ij}$ avec $1\leq i\leq \ell$ et $1\leq j\leq n$. \'{E}tant donn\'ee une famille $\frak{S}_n$-stable $F$ des 
polyn\^{o}mes homog\`{e}nes en les variables $x_{ij}$, le plus petit espace vectoriel $\mathcal{M}_F$ clos par d\'erivation partielle et clos par l'action des op\'erateurs de polarisation contenant $F$ est {le module de polarisation} engendr\'e par la famille $F$. Les modules $\mathcal{M}_F$ sont tous des repr\'esentations du produit direct $\frak{S}_n\times GL_{\ell}(\mathbb{C})$. Dans le but d'\'etudier la d\'ecomposition en sous-modules irr\'eductibles on calcule la caract\'eristique de Frobenius gradu\'ee de ces modules. Pour plusieurs cas de familles homog\`{e}nes $\mathfrak{S}_n$-stables constitu\'ees des polyn\^{o}mes homog\`{e}nes \`{a} $n$ variables, pour tout $n\geq 1$, on d\'emontre des formules g\'en\'erales pour cette caract\'eristique gradu\'ee de fa\c con globale, ind\'ependante de la valeur de $\ell$.

\end{abstract}

\section{Introduction}\label{Intro}

This work is inspired by the Operator Theorem (ex-operator conjecture) of M. Haiman (see \cite{Haiman3}). This theorem states that the smallest subspace of $\complexes[x_1,\ldots,x_n]$ closed under taking partial derivatives ${\frac{\partial }{\partial x_i}}$, closed under the action of generalized polarization operators, $E_{p}=\sum_{j=1}^{n}y_{j}\frac{\partial^p }{\partial x_j^p}$, that contains the Vandermonde determinant $\Delta_{n}({\bold x})~:=~\prod_{1\leq i<j\leq n}(x_{i}-x_{j})$, coincides with the space $\mathcal{D}_n$ of diagonal harmonics polynomials of $\mathfrak{S}_n$. The space $\mathcal{D}_n$ consist of all polynomials $f$ in the variables $x_1,\ldots,x_n,y_1,\ldots,y_n$ killed by the power sum differential operators $\sum_{j=1}^{n}\frac{\partial^{h} }{\partial x_j^{h}}\frac{\partial^{k} }{\partial y_j^{k}}$ with $1\leq h+k\leq n$. In others words this space is generated by the Vandermonde determinant as a module over the algebra of operators $\complexes\left[\frac{\partial }{\partial x_1},\ldots,\frac{\partial }{\partial x_n},E_1,\ldots,E_{n-1}\right]$. In this work we generalize this construction by considering polynomials in $\ell$ sets of $n$ variables, that is, polynomials in the variables $x_{11},\ldots,x_{1n},\ldots,x_{\ell 1},\ldots,x_{\ell n}$. We are interested in the decomposition into irreducible submodules of these spaces. 

We start by generalizing the context of the Operator Theorem to the context of polynomials in the matrix variables $X=(x_{ij})$, with $1\leq i\leq \ell$ and $1\leq j\leq n$. The diagonal action of $\mathfrak{S}_n$ on this polynomials in $\ell$ sets of $n$ variables is defined by permuting the columns of $X$, that is, for any permutation $\sigma\in\mathfrak{s}_n$ we replace the variable $x_{ij}$ by $x_{i\sigma(j)}$. We say that a family $F$ of homogeneous polynomials (in $X=(x_{ij})$) is $\mathfrak{S}_n$-stable if $F$ is closed under the diagonal action of $\mathfrak{S}_n$. Given any such a family $F$, we define the polarization module generated by the family $F$ as the smallest vector space closed under taking partial derivatives $\frac{\partial }{\partial x_{ij}}$, closed under the action of generalized polarization operators $E_{i,k}^{(p)}=\sum_{j=1}^{n}x_{ij}\frac{\partial^{p} }{\partial x_{kj}^{p}}$, that contains $F$. The diagonal action of $\mathfrak{S}_n$ makes $\mathcal{M}_F$ an $\mathfrak{S}_n$-module. The closure by the action of polarization operators $E_{i,k}^{(p)}$ is equivalent to the closure by the action $x_{ij}\longmapsto \sum_{k=1}^{\ell}m_{ik}x_{kj}$ where $M=(m_{ij})\in{GL}_{\ell}(\complexes)$ (see \cite{ProcesiKraft,Procesi}). Then, with this action $\mathcal{M}_F$ is also a polynomial representation of ${GL}_{\ell}(\complexes)$. The two actions of $\mathfrak{S}_n$ and ${GL}_{\ell}(\complexes)$ commutes and this implies that $\mathcal{M}_{F}$ is a representation of the direct product $\mathfrak{S}_n\times{GL}_{\ell}(\complexes)$. In particular, when the family $F$ is the orbit of a single homogeneous polynomial $f$, that is, $F=\{\sigma\cdot f \ \vert \ \sigma\in\mathfrak{S}_n\}$ we denote the polarization module generated by $F$ simply as $\mathcal{M}_f$. Also, we call $\mathcal{M}_f$ the polarization module generated by $f$. 

Particular cases of this construction correspond to certain important spaces in Combinatorics (see \cite{FBergeron}) and Algebraic Geometry (see \cite{Geramita}). For instance, when $\ell=1$ (1 set of variables ${\bold x}=x_1,\ldots,x_n$) and $f=\Delta_{n}({\bold x})$ we get the space $\mathcal{H}_{n}$ of harmonics polynomials of the group $\mathfrak{S}_n$, that is, the space of polynomials zeros of the power sum differential operators $\sum_{i=1}^{n}\frac{\partial}{\partial x_{i}^{k}}$ with $k$ such that $1\leq k\leq n$, (see \cite{FBergeron2}). It's remarkable that the dimension of the space $\mathcal{H}_n$ is $n!$ (see \cite{Humphreys}). As we already said, the case $\ell=2$ and $f=\Delta_{n}({\bold x})$ gives us the space $\mathcal{D}_n$ of diagonal harmonics polynomials. It's also a remarkable fact that $\dim(\mathcal{D}_n)=(n+1)^{n-1}$ (the number of parking functions) (see \cite{Haiman3}). The space of {trivariate diagonal harmonics $\mathcal{D}_{n}^{(3)}$} consist of polynomials zeros of the polarized power sums operators 
$\sum_{j=1}^{n}\frac{\partial^{a}}{\partial {x_j}^{a}}\frac{\partial^{b}}{\partial {y_j}^{b}}\frac{\partial^{c} }{\partial {z_j}^{c}}$ (see \cite{FBergeronEminem}), where $1\leq a+b+c\leq n$. M. Haiman conjectured that this space coincides with the polarization module generated by the Vandermonde determinant when $\ell=3$ (see \cite{Haiman1}) this is still an open problem. Also, M. Haiman conjectured that $\dim(\mathcal{D}_{n}^{(3)})=2^{n}(n+1)^{n-2}$. F. Bergeron extended this conjecture to the space $\mathcal{H}_{n}^{(r)}$ of higher diagonal harmonics polynomials which has dimension 
$\dim(\mathcal{H}_{n}^{(r)})=(r+1)^{n}(rn+1)^{n-2}$ (\cite{FBergeronEminem}). Also, F. Bergeron extend Haiman's conjecture concerning to the identification of the space $\mathcal{D}_n^{(\ell)}$ of multivariate diagonal harmonics polynomials with polarization module generated by the Vandermonde determinant for $\ell>3$.  

The goal of this paper is to study the decomposition into irreducible submodules under the action of $\mathfrak{S}_n\times{GL}_{\ell}(\complexes)$ of polarization modules $\mathcal{M}_F$. To do this we compute explicitly the graded Frobenius characteristic of $\mathcal{M}_F$ in the form
\begin{equation}\label{Formula1Introduccion}
\mathcal{M}_{F}({\bold q},{\bold w})=\sum_{\lambda\vdash n}
\,\sum_{\vert\mu\vert\leq d}b_{\lambda,\mu}s_{\mu}({\bold q})s_{\lambda}({\bold w}) 
\end{equation}
where $b_{\lambda,\mu}\in\naturals$, ${\bold q}=q_1,q_2,\ldots,q_{\ell}$, ${\bold w}=w_1,w_2,\ldots$, and $d$ is the maximal degree of polynomials in $F$. Here the Schur functions $s_{\mu}({\bold q})$ encodes irreducible for $GL_{\ell}(\complexes)$ while $s_{\lambda}({\bold w})$ encodes irreducible for $\mathfrak{S}_n$. The coefficients $b_{\lambda,\mu}$ are the multiplicities of irreducible submodules under the action of $\mathfrak{S}_n\times{GL}_{\ell}(\complexes)$. A theorem of F. Bergeron shows that this coefficients are independent of $\ell$, in others words, the value of $\ell$ only affects the formula (\ref{Formula1Introduccion}) giving the numbers of ${\bold q}$ variables appearing in $s_{\mu}({\bold q})$ (see \cite{FBergeron}). Furthermore, he shows in \cite{FBergeron} that $\mu$ has at most $n$ parts. This leads us to obtain a general formula for (\ref{Formula1Introduccion}) that holds for any  $\ell\geq 1$ if we can compute it for every $\ell\leq n$. 

In this paper we completely describe the decomposition into irreducible submodules of the polarization modules generated by each of the polynomials $p_1^d$, $p_d$, and $e_d$ for any $d\geq 1$. We construct an explicit linear basis of each module and then we compute the graded Frobenius characteristic of $\mathcal{M}_{p_1^d}$, $\mathcal{M}_{p_d}$ and $\mathcal{M}_{e_d}$. We propose a conjectural formula for the graded Frobenius series of the polarization module generated by the monomial symmetric function $m_{(2,1^{d-2})}$, in any degree $d$. For instance, we also conjecture that, for $d\geq 5$ we must have
$\mathcal{M}_{e_{d-1,1}}\cong\mathcal{M}_{m_{(2,1^{d-2})}}$ as $\mathfrak{S}_n\times{GL}_{\ell}(\complexes)$-modules. We believe that these are modules are the elementary building blocks for a general classification of polarization modules generated by a given homogeneous symmetric polynomial in any number of variables $n$. Experimental evidence has shown that the last assertion is true up to degree 5 and $n\leq 6$. In particular, we completely determine the classification of polarization modules generated by a single homogeneous symmetric polynomial when the degree is 2 or 3. For the case of degree 4 and 5 we have a conjectural classification that seems to be complete. This framework lead us to think that the Hilbert series of the polarization modules generated by a single homogeneous symmetric polynomial, in any degree, is always $h$-positive (see Conjecture \ref{Hpositividad}). 

Obviously, we have $\mathcal{M}_f\cong\mathcal{M}_{k\cdot f}$ for every scalar $k$. In order to classify up to isomorphism polarization modules generated by a given homogeneous symmetric polynomial of any degree $d$, we identify any non zero homogeneous symmetric polynomial $f$ of degree $d$, written in the monomial basis as $f~=~\sum_{\lambda\vdash d}c_{\lambda}m_{\lambda}$, with a point in the real projective space $\mathbb{RP}^{p(d)-1}$, where $p(d)$ is the number of integer partitions of $d$. The homogeneous coordinates of the corresponding  point are ordered according to the following order on integer partitions of $d$:
$(d)$, $(d-1,1)$, $(d-2,2)$, $(d-2,1,1)$,$\ldots$, $(1,1,\ldots,1)$. In degree 2, we show that there are two types of polarization modules 
$\mathcal{M}_{p_1^2}$ and $\mathcal{M}_{p_2}$ up to isomorphism. More precisely, if $[a:b]\in\mathbb{RP}^1$ and $f=a\cdot m_{2}+b\cdot m_{11}$ then $\mathcal{M}_{f}\cong\mathcal{M}_{p_1^2}$ when $[a:b]=[1:2]$, while $[a:b]\neq [1:2]$ implies $\mathcal{M}_{f}\cong\mathcal{M}_{p_2}$.
Notice that the last statement is independent on the number of variables $n$. The situation when the degree is 3 is more complicated, in this case, we will need to introduce the notion of $n$-exception to completely classify these polarization modules. Let $n\geq 3$, a point $[a:b:c]\in\mathbb{RP}^{2}$ is a $n$-exception if and only if\ $[a:b:c]\neq[1:3:6]$ and $6a(2b+(n-2)c)=4(n-1)b^2$. When $n=2$, $[a:b:c]$ is a 2-exception if and only if $b=0$ or $b=3a$. There are three types of polarization modules generated by a single polynomial of the form $f=a\cdot m_3+b\cdot m_{21}+c\cdot m_{111}$. If $[a:b:c]=[1:3:6]$ then $\mathcal{M}_f\cong\mathcal{M}_{p_1^3}$; if $[a:b:c]$ is a $n$-exception then $\mathcal{M}_f\cong\mathcal{M}_{p_3}$; otherwise, 
$\mathcal{M}_f\cong\mathcal{M}_{h_3}$. These results are valid for any $\ell$ (the number of sets of $n$ variables). Also, we will see that $n$-exceptions appear in any degree $d\geq 3$.  Characterize $n$-exceptions for degrees higher than 3 is a problem for the future. 

For particulars $\mathfrak{S}_n$-stable families of homogeneous polynomials we completely describe the graded Frobenius characteristic of the associated polarization module. In particular, we have a conjecture giving constraints for the multiplicities of irreducible submodules of polarization modules generated by any $\mathfrak{S}_n$-stable family consisting of homogeneous polynomials of degree at most 2. Similarly, for any $\mathfrak{S}_n$-stable families of homogeneous polynomials of degree at most 3. More precisely, we compute the graded Frobenius characteristic of polarization modules generated by the family of all monomials of degree $d$, when $d=2,3$ in any number of variables $n$. Notice that all these formulas holds for any $\ell$.

\section{Preliminaries}
Let $X:=(x_{ij})_{i,j}$ be a $\ell\times n$ matrix of commuting and independent variables $x_{ij}$
For any $i$, we call the $i^{th}$-row of $X$, denoted by ${\bold x}_{i}:=(x_{i1},\ldots,x_{in})$, the $i^{th}$ set of variables. For any $j$, $X_{j}$ denotes the $j^{th}$ column of $X$. The same convention is adopted for any $\ell\times n$ non negative integer matrix of exponents $A$. Then the monomials are defined as follows:
\begin{equation}
X^{A}:=x_{11}^{a_{11}}\cdots x_{ij}^{a_{ij}}\cdots x_{\ell{n}}^{a_{\ell{n}}}.
\end{equation}
these monomials form a linear basis of the $\complexes$-vector space $\mathcal{R}_{n}^{(\ell)}=\complexes[X]$ of polynomials in $\ell$ sets of $n$ variables. The (vector) \textbf{degree} $\deg\big(X^{A}\big)$ lies in $\mathbb{N}^{\ell}$ and is given by
$\deg\left(X^{A}\right):=\left(\sum_{j=1}^{n}a_{1j},\dots,\sum_{j=1}^{n}a_{\ell j}\right)$. For each ${\bold d}\in\mathbb{N}^{\ell}$, we denote by $\mathcal{R}_{n,{\mathbf d}}^{(\ell)}$ the span of degree ${\mathbf d}$ monomials in $\mathcal{R}_{n}^{(\ell)}$. Then $\mathcal{R}_{n}^{(\ell)}$ is a $\mathbb{N}^{\ell}$-graded, that is, $\mathcal{R}_{n}^{(\ell)}=\bigoplus_{{\mathbf d}\in\naturals^{\ell}}\mathcal{R}_{n,{\mathbf d}}^{(\ell)}$. We will consider homogeneous subspace $\mathcal{V}$ of $\mathcal{R}_{n}^{(\ell)}$, the degree ${\bold d}$ homogeneous component of $\mathcal{V}$ is denoted by $\mathcal{V}_{\bold d}$. Recall that, $\mathcal{V}_{\bold d}:=\mathcal{V}\cap\mathcal{R}_{n,{\bold d}}^{(\ell)}$ and the \textbf{Hilbert series} of $\mathcal{V}$ is $\mathcal{V}({\bold q}):=\sum_{{\bold d}\in\mathbb{N}^{\ell}}\dim(\mathcal{V}_{\bold d}){\bold q}^{\bold d}$, where ${\bold q}:=(q_1,\ldots,q_{\ell})$ and ${\bold q}^{\bold d}:=q_{1}^{d_1}\cdots q_{\ell}^{d_{\ell}}$. \\

On $\mathcal{V}$ we consider two linear group actions:
\begin{enumerate}
\item The (left) diagonal action of $\mathfrak{S}_n$, \ \ $\sigma\cdot X^A:=x_{1\sigma(1)}^{a_{11}}\cdots x_{i\sigma(j)}^{a_{ij}} \cdots
x_{{\ell}\sigma(n)}^{a_{{\ell}n}}$,\ $\forall\sigma\in\frak{S}_n$,
\item The (right) action of $GL_{\ell}(\mathbb{C})$, \ \ $X^A\cdot M :=(MX)^A$,\ $\forall M\in GL_{\ell}(\mathbb{C})$, that is, $x_{ij}\longmapsto \sum_{k=1}^{\ell}m_{ik}x_{kj}$, for every matrix $M=(m_{ij})\in GL_{\ell}(\mathbb{C})$.
\end{enumerate}
These two group actions on $\mathcal{V}$ commute and then we can consider $\mathcal{V}$ as a representation of the direct product $\mathfrak{S}_n\times{GL}_{\ell}(\complexes)$, with the (left) action, $(\sigma,M)\cdot X^{A}:=\sigma\cdot (M^{-1}X)^{A}$. Then we have a direct sum decomposition of the form (see \cite{BruceSagan,FultonHarris,Procesi}) 
\begin{equation}\label{DecompIrred}
\mathcal{V}=\bigoplus_{\lambda\vdash n}\bigoplus_{\mu} b_{\lambda,\mu}\,\mathcal{W}_{\mu} \otimes \mathcal{U}_{\lambda},
\end{equation}
where $b_{\lambda,\mu}\in\naturals$, the\ $\mathcal{U}_{\lambda}$ are irreducible $\frak{S}_n$-modules and the $\mathcal{W}_{\mu}$ are irreducible polynomial representations of  $GL_{\ell}(\mathbb{C})$. The \textbf{graded Frobenius characteristic} of $V$ is defined as follows
\begin{equation}
\mathcal{V}({\bold q},{\bold w}):=\sum_{{\bold d}\in\mathbb{N}^{\ell}}\left(
\frac{1}{n!}\sum_{\sigma\in\frak{S}_n}{\chisotazo}_{\mathcal{V}_{\bold d}}(\sigma)\,
\,p_{_{\lambda(\sigma)}}({\bold w})\right){\bold q}^{\bold d},
\end{equation}
where $\chisotazo_{\mathcal{V}_{\bold d}}$ is the $\mathfrak{S}_n$-character of $\mathcal{V}_{\bold d}$ and $p_{\lambda(\sigma)}({\bold w}):=p_{1}({\bold w})^{c_1(\sigma)}\cdots{p_{n}({\bold w})^{c_{n}(\sigma)}}$. It's well known (see \cite{FBergeron}) that the graded Frobenius characteristic of $\mathcal{V}$ 
has the following form (not depending on $\ell$):
\begin{equation}\label{FormaGeneralDelaCaracteristicaDeFrobenius}
\mathcal{V}({\bold q},{\bold w})=\sum_{\lambda\vdash n}\sum_{\mu}b_{\lambda,\mu}s_{\mu}({\bold q})s_{\lambda}({\bold w}),
\end{equation}
where $b_{\lambda,\mu}$ are the multiplicities in (\ref{DecompIrred}). $\ell(\mu)\leq n$ (see \cite{FBergeron} for more details). The Schur functions $s_{\mu}({\bold q})$ encode the irreducible polynomial representations of $GL_{\ell}(\complexes)$ in $\mathcal{V}$ of type $\mu$, and the Schur functions $S_{\lambda}({\bold w})$ encodes the irreducible $\mathfrak{S}_n$-modules in $V$ of type
$\lambda$. Recall that the Hilbert Series of $\mathcal{V}$ is obtained by replacing $s_{\lambda}({\bold w})$ by $f^{\lambda}$ (the Hook length formula) for every $\lambda$ in formula \ref{FormaGeneralDelaCaracteristicaDeFrobenius}.


\section{Definitions and discussions}\label{Definitions}
We denote the partial derivative operator on $\mathcal{R}_{n}^{(\ell)}$ by $\displaystyle{\partial_{ij}:=\frac{\partial \ \ }{\partial x_{ij}}}$. We use the \textbf{generalized polarization operators $E_{i,k}^{(p)}:\mathcal{R}_n^{(\ell)}\longrightarrow\mathcal{R}_{n}^{(\ell)}$} given by $\displaystyle{E_{i,k}^{(p)}:=\sum_{j=1}^{n}x_{ij}\partial_{kj}^{p}}$
where $1\leq i,k\leq \ell$ and $p\geq 1$. Clearly the operators $E_{i,k}^{(p)}$ depend on the choice of $(i,k,p)$. For $p=1$ we simply write $E_{i,k}:=E_{i,k}^{(1)}$ \  (see, \cite{HWeyl,Hunziker}, and \cite{Procesi} for more details). 
\subsection{Generalized Polarization Modules}
We say that a subspace $\mathcal{V}$ of $\mathcal{R}_{n}^{(\ell)}$ is \textbf{closed under derivatives} if for every $g\in\mathcal{V}$ we have $\partial_{ij}(g)\in\mathcal{V}$, for all $(i,j)$ such that $1\leq i\leq \ell$ and $1\leq j\leq n$. We say that $\mathcal{V}$ is \textbf{closed under polarization} if $E_{i,k}^{(p)}(g)\in \mathcal{V}$ for all $g\in \mathcal{V}$ and all suitable triple $(i,k,p)$. Since any intersection of closed under derivatives subspaces is closed under derivatives we define the \textbf{derivative closure $\boldsymbol{\mathcal{D}}(\mathcal{V})$ of $\mathcal{V}$} as the smallest subspace of $\mathcal{R}_{n}^{(\ell)}$ closed under derivatives that contains $\mathcal{V}$. Similarly, we define the \textbf{polarization closure $\boldsymbol{\mathcal{E}}(\mathcal{V})$ of $\mathcal{V}$}. It's not hard to show the following result (proofs will be avaible in a complete version of this paper, see \cite{Blandin}):
\begin{lemma}
Let $\mathcal{V}$ be a homogeneous subspace of $\mathcal{R}_n^{(\ell)}$ then $\boldsymbol{\mathcal{E}}(\boldsymbol{\mathcal{D}}(\mathcal{V}))=\boldsymbol{\mathcal{D}}(\boldsymbol{\mathcal{E}}(\mathcal{V}))$.
\end{lemma}
We set $\boldsymbol{\mathcal{P}}(\mathcal{V}):=\boldsymbol{\mathcal{E}}(\boldsymbol{\mathcal{D}}(\mathcal{V}))$.
\noindent A subset $F$ of $\mathcal{R}_{n}^{(\ell)}$ is called a \textbf{homogeneous stable family} if the following conditions holds:
\begin{enumerate}
\item $F$ consist only of homogeneous polynomials, 
\item $F$ is \textbf{stable} (or $\mathfrak{S}_n$-stable w.r.t. the diagonal action of $\mathfrak{S}_n$), that is, for any permutation $\sigma\in\mathfrak{S}_n$ we have $\sigma\cdot g\in F$, for all $g\in F$.
\end{enumerate}
\begin{defi}
For a given homogeneous stable family $F$,  we set $\mathcal{M}_{F}$ to be the smallest $\mathbb{R}$-vector space closed under derivatives and closed under polarization containing the family $F$. We call the vector space $\mathcal{M}_{F}$ the \textbf{polarization module generated by the family $F$}.
\end{defi}
We can describe the polarization module generated by $F$ as
$\mathcal{M}_{F}:=\boldsymbol{\mathcal{P}}\big(\mathbb{R}\cdot F\big)$,
where $\mathbb{R}\cdot F$ denotes the real vector space spanned by $F$.
\begin{remark}
When the family $F$ consist of only one homogenous polynomial $f\in\mathcal{R}_{n}^{(\ell)}$, we denote by $\mathcal{M}_f$ the polarization module generated by the family $\{\sigma\cdot f\ \vert\ \sigma\in\mathfrak{S}_n\}$. In symbols,
\[ \mathcal{M}_{f}:=\boldsymbol{\mathcal{P}}\big(\mathbb{R}\cdot\{\sigma\cdot f\ \vert\ \sigma\in\mathfrak{S}_n\}\big).\]
\end{remark}

\subsection{Properties of Polarization Modules}
The closure under polarization operators property implies that any polarization modules is a polynomial representation of $GL_{\ell}(\complexes)$ with the (right) action $f(X)\cdot M:=f(MX)$. In fact 
\begin{lemma}[see C. Procesi, \cite{Procesi}]
A subspace $V$ of $\mathcal{R}_{n}^{(\ell)}$ is closed under the action of polarization operators $E_{i,k}$ (when $p=1$) if and only if $V$ is a $GL_{\ell}(\mathbb{C})$-module with the action $X^{A}\cdot M=(MX)^{A}$.
\end{lemma}
For any $g\in\mathcal{R}_{n}^{(\ell)}$ we have the identities:\ 
$\sigma\cdot E_{i,k}^{(p)}(g)=E_{i,k}^{(p)}(\sigma\cdot g)$ and  
$\sigma\cdot \partial_{ij}(g)=\partial_{i,\sigma(j)}(\sigma\cdot g)$. These identities implies the following assertions for any family $F$
\begin{lemma}
$\mathcal{M}_{F}$ is a representation of $\frak{S}_n$ with the diagonal action of $\mathfrak{S}_n$.
\end{lemma}
Recall that the two actions of $\mathfrak{S}_n$ and $GL_{\ell}(\mathbb{C})$ on $\mathcal{R}_{n}^{(\ell)}$ commutes and so, we can assert that
\begin{lemma}
$\mathcal{M}_{F}$ is a\ $\frak{S}_{n}\times GL_{\ell}(\mathbb{C})$-module with the action $(\sigma,M)\cdot f(X):=f(M^{-1}X_{\sigma})$,
where $X_{\sigma}:=(x_{i\sigma(j)})_{i,j}$.
\end{lemma}
\section{Frobenius characteristics of some polarization modules}
In the following lines we use the notation of Macdolnald's book \cite{IGMacdonald}. The proofs of the following results will be available in a complete version of this paper (see \cite{Blandin}).
\begin{teor}\label{ProposicionUNO}
Let $d$ be a positive integer. The following formulas holds for any $\ell\geq 1$
\begin{equation}\label{ProposicionUNOPrimeraParte}
\mathcal{M}_{e_{1}^{d}}({\bold q},{\bold w})=\left(\sum_{j=0}^{d}s_{j}({\bold q})\right)s_{n}({\bold w})=\left(\sum_{j=0}^{d}h_{j}({\bold q})\right)h_{n}({\bold w}).
\end{equation}
\begin{align}\label{ProposicionUNOSegundaParte}
\mathcal{M}_{p_{d}}({\bold q},{\bold w})&=\left(\sum_{j=0}^{m}s_{j}({\bold q})\right)\cdot {S}_{n}({\bold w})
+\left(\sum_{j=1}^{m-1}s_{j}({\bold q})\right)\cdot {S}_{n-1,1}({\bold w})\\
&=\big(1+h_{d}({\bold q})\,\big)\,h_{n}({\bold w})+\left(\sum_{j=1}^{d-1}h_{j}({\bold q})\right)h_{n-1,1}({\bold w}).
\end{align}
\begin{align}\label{ProposicionUNOTerceraParte}
\mathcal{M}_{e_{d}}({\bold q},{\bold w})
&=\sum_{i=0}^{\lfloor{d/2}\rfloor}\left(\sum_{j=i}^{d-i}s_{j}({\bold q})\right)S_{n-i,i}({\bold w})\\
&=\sum_{i=0}^{\lfloor\frac{d}{2}\rfloor}h_{n-i,i}({\bold w})\,h_{i}({\bf q})
+\sum_{i=\lfloor\frac{d}{2}\rfloor+1}^{d}h_{n-d+i,d-i}({\bold w})\,h_{i}({\bf q}).
\end{align}
\end{teor}
\begin{teor}\label{ProposicionDOS}
Let $f({\bold x}_1)$ be a symmetric polynomial of degree 2  in $n\geq 2$ variables ${\bold x}_1=x_{11},x_{12},\ldots,x_{1n}$. Suppose that $f({\bold x}_1)$ is given in the monomial basis as follows:
\begin{equation}
f({\bold x}_1)=a\cdot m_{2}({\bold x}_1)+b\cdot m_{1,1}({\bold x}_1),
\end{equation}
then the Frobenius characteristic of the space $\mathcal{M}_{f}$ is given by one of the following two cases:
\begin{equation}
\mathcal{M}_{f}({\bold{q}},{\bold{w}})=
\left\{
\begin{array}{cc}
\big(1+h_1({\bold q})+h_2({\bold q})\big)\cdot h_{n}({\bold{w}}) &  \text{if} \  \ \big[a:b\big]=\big[1:2\big], \\
\big(1+s_1({\bold q})+s_2({\bold q})\big)\cdot s_{n}({\bold{w}})+s_1({\bold q})\cdot s_{n-1,1}({\bold{w}})\\
=\big(1+h_{2}({\bold{q}})\,\big)\cdot h_{n}({\bold{w}})+h_{1}({\bold q})\cdot h_{n-1,1}({\bold w})
& {\rm otherwise}.
\end{array}
\right.
\end{equation}
\end{teor}
\begin{corollary}\label{CorolarioGrado2}
If $f$ is a homogeneous symmetric polynomial of degree 2, then the associated $\frak{S}_n$-module $\mathcal{M}_f$ is isomorphic as an $\mathfrak{S}_n\times GL_{\ell}(\mathbb{K})$-module to one of the two modules $\mathcal{M}_{p_1^{2}}$,
$\mathcal{M}_{p_{2}}$.
\end{corollary}
\subsection{Exceptions}
Recall that, we identify $f=a\cdot m_3+b\cdot m_{21}+c\cdot m_{111}$, ($a$,$b$,$c$ in $\reals$) with it's homogeneous coordinates 
$[f]:=[a:b:c]\in\mathbb{RP}^{2}$. For instance $p_{1}^{3}$ is the point $[1:3:6]$, $p_{21}$ is $[1:1:0]$ and $h_{3}$ is the point $[1:1:1]$. 
\begin{defi}
We say that an homogeneous symmetric polynomial $f$ in 
$\reals[x_{11},\ldots,x_{1n}]$ is a \textbf{$n$-exception} if 
\[ \dim\left(\reals\{\partial_{11}(f),\ldots,\partial_{1n}(f),E_{1,1}^{(2)}(f)\}\right)=n. \]
In others words, $f$ is a $n$-exception if the dimension of the real linear span of its first order partial derivatives $\partial_{11}(f),\ldots,\partial_{1n}(f)$ and the polynomial $E_{1,1}^{(2)}(f)$ has dimension $n$.  
\end{defi}
For example, the points $[1:0:0]$ ( $f=p_3$ ) and  $[0:0:1]$ ( $f=e_3$) are $n$-exceptions for every $n\geq 2$. For instance, $[3:3:-2]$ is a 3-exception, $[9:21:28]$ is a $4$-exception, $[2:3:2]$ is a $5$-exception, $[4:-3:4]$ is a $5$-exception. Another example is $[1,1,0]$ which is a $4$-exception, because $E_{1,1}^{(2)}p_{21}=\sum_{j=1}^{4}\partial_{{1j}}p_{21}$. 
For instance, for every $n\geq 3$, the point $[1:1:1]$ ($f=h_{3}$) is not an $n$-exception, the point $[0:1:0]$ ($f=m_{21}$) is not an $n$-exception. En degree 4, the point $[5:14:21:28:35]$ is a 11-exception.
\begin{teor}
Let $f$ be a homogeneous symmetric polynomial of degree 3 in $\reals[x_{11},\ldots,x_{1n}]$. If $[f]\neq[1:3:6]$ then 
\[ \dim\left(\reals\{\partial_{11}(f),\ldots,\partial_{1n}(f),E_{1,1}^{(2)}(f)\}\right)\geq n. \]
\end{teor}
\begin{teor}
Suppose that $n\geq 3$ and $f({\bold x}_1)=a\cdot m_3({\bold x}_1)+b\cdot m_{21}({\bold x}_1)+c\cdot m_{111}({\bold x}_1)$, ($a$,$b$,$c$ dans $\reals$). Then $f$ is a $n$-exception if and only if $[a:b:c]\neq[1:3:6]$ and $6a(2b+(n-2)c)=4(n-1)b^2$. When $n=2$, $[a:b:c]$ is a 2-exception if and only if $b=0$ or $b=3a$. 
\end{teor}
\begin{teor}\label{ProposicionTRES}
Let $f$ be a homogeneous symmetric polynomial of degree 3 in $n\geq 2$ variables. Suppose that $f\in\reals[x_{11},\ldots,x_{1n}]$ and $[f]=[a:b:c]$, then the Frobenius characteristic of the $\frak{S}_n$-module $\mathcal{M}_{f}$ is given by one of the following three cases:
\begin{align}
&\begin{cases}
\big(\,1+h_1({\bold q})+h_2({\bold q})+h_3({\bold q})\,\big)\cdot h_{n}({\bold w}) & \text{if} \ \ \big[a:b:c\big]=\big[1:3:6\big],\\
\big(\,1+h_3({\bold q})\,\big)\cdot h_{n}({\bold{w}})+\big(\,h_{1}({\bold q})+h_{2}({\bold q})\,\big)h_{n-1,1}({\bold{w}}) &\\
=(1+s_1({\bold q})+s_2({\bold q})+s_3({\bold q}))\cdot s_{n}({\bold w})
+(s_1({\bold q})+s_2({\bold q}))\cdot s_{n-1,1}({\bold w})& \text{if} \ \ \big[a:b:c\big]\ \text{is an $n$-exception,}\\
\big(\,1+h_{2}({\bold q})+h_3({\bold q})\,\big)\cdot h_{n}({\bold{w}})+\big(h_1({\bold q})+h_2({\bold q})\,\big)\cdot h_{n-1,1}({\bold{w}}) & \\
=(1+s_1({\bold q})+\boldsymbol{2}s_2({\bold q})+s_3({\bold q}))\cdot s_{n}({\bold w})+(s_1({\bold q})+s_2({\bold q}))\cdot s_{n-1,1}({\bold w})&\ \ \text{otherwise}.
\end{cases}
\end{align}
\end{teor}
\vspace{0.35cm}
\begin{corollary}\label{CorolarioGrado3}
Le $f$ be a homogeneous symmetric polynomial of degree 3. There are three types of polarization modules $\mathcal{M}_{p_1^{3}}$, $\mathcal{M}_{p_{3}}$ \ or $\mathcal{M}_{h_3}$.
\end{corollary}
\section{Conjectures}
\begin{conj}
For any degree $d\geq 2$, $f=p_{2}p_{1}^{d-2}$ is a $(d+1)$-exception.
\end{conj}
\begin{conj}
Let $f$ be a homogeneous symmetric polynomial of degree $d$ in $\mathbb{R}[x_{11},\ldots,x_{1n}]$. If $f$ is not a scalar multiple of $p_{1}^{d}$ then
\begin{equation*}
\dim\left(\reals\{\partial_{11}(f),\ldots,\partial_{1n}(f),E_{1,1}^{(2)}(f)\}\right)\geq n.
\end{equation*}
\end{conj}
In the following, we will write $s_{j}:=s_{j}({\bold q})$ for short.
\begin{conj}\label{ConjeturaM2111}
Let $d\geq 3$ be an integer. Suppose that $m_{_{211\cdots1}}$ is the monomial symmetric function indexed by the partition of $d$ $\mu=(2,1^{d-2})$. If $n\geq d$ we have 
\begin{align*}
&\mathcal{M}_{m_{_{211\cdots1}}}({\bold q},{\bold w})
=\left(1+s_1+\boldsymbol{2}\cdot\sum_{j=2}^{d-1}s_{j}+s_d\right)s_{n}({\bold w})
+\sum_{i=1}^{\lfloor\frac{d}{2}\rfloor-1}\left(s_i+\boldsymbol{2}\cdot\sum_{j=i+1}^{d-i-1}s_j
+s_{d-i}\right)s_{n-i,i}({\bold w})\\
&+
\left(\sum_{j=\lfloor\frac{d}{2}\rfloor}^{d-\lfloor\frac{d}{2}\rfloor}s_{j}\right)s_{n-\lfloor\frac{d}{2}\rfloor,\lfloor\frac{d}{2}\rfloor}({\bold w}).
\end{align*}
Also, if $d\geq 5$, we believe that the following isomorphism as a representation of $\mathfrak{S}_n\times{GL}_{\ell}(\complexes)$ holds
\[ \mathcal{M}_{m_{2,1,1,1,\cdots 1}}\cong\mathcal{M}_{e_{d-1,1}}.\]
\end{conj}
Consider the following $\frak{S}_n$-stable families of homogeneous polynomials in the variables $x_{11},\ldots,x_{1n}$ $\mathcal{A}:=\big\{x_{1j}^d\big\vert\ 1\leq j\leq n\big\}$,\ $\mathcal{B}:=\big\{x_{1,i}^d-x_{1,j}^d\big\vert\ 1\leq i<j\leq n \big\}$ \ and \ $\mathcal{C}:=\left\{\,\prod_{a\in A}x_{1,a}\,\big\vert\ A\subseteq [n],\ \vert A\vert=d\right\}$.
\begin{conj}\label{FamiliasABC}
The graded Frobenius characteristic of the families $\mathcal{A}$, $\mathcal{B}$ and $\mathcal{C}$ are given by 
\begin{equation*}
\mathcal{M}_{_\mathcal{A}}({\bold q},{\bold w})
=\left(\sum_{j=0}^{d}s_{j}({\bold q})\right)s_{n}({\bold w})
+\left(\sum_{j=1}^{d}s_{j}({\bold q})\right)s_{n-1,1}({\bold w}).
\end{equation*}
\begin{equation*}
\mathcal{M}_{_\mathcal{B}}({\bold q},{\bold w})
=\left(\sum_{j=0}^{d-1}s_{j}({\bold q})\right)s_{n}({\bold w})
+\left(\sum_{j=1}^{d}s_{j}({\bold q})\right)s_{n-1,1}({\bold w}).
\end{equation*}
\begin{equation*}
\mathcal{M}_{_\mathcal{C}}({\bold q},{\bold w})=\left(\sum_{j=0}^{d}s_{j}({\bold q})\right)s_{n}({\bold w})
+\sum_{i=1}^{\lfloor{d/2}\rfloor}\left(\sum_{j=i}^{d-i+1}s_{j}({\bold q})\right)s_{n-i,i}({\bold w}).
\end{equation*} 
\end{conj}
\begin{conj}\label{ConjeturaGrado2Familia}
Let $\mathcal{T}_2$ be the family of all monomials of total degree 2 in the variables $x_{11},\ldots,x_{1n}$. Then the graded Frobenius characteristic of its associated polarization module is 
\begin{equation*}
\mathcal{M}_{\mathcal{T}_2}({\bold q},{\bold w},n)
=\begin{cases}
(1+s_1+s_2)s_{1}({\bold w}) & \text{if}\ n=1,\\
(1+s_1+\boldsymbol{2}\,s_2)s_{2}({\bold w})+(s_1+s_2)s_{1,1}({\bold w}) & \text{if}\ n=2, \\
(1+s_1+\boldsymbol{2}\,s_2)s_{3}({\bold w})+(s_1+\boldsymbol{2}\,s_2)s_{2,1}({\bold w}) &  \text{if}\ n=3,\\
\left(1+s_{{1}}+\boldsymbol{2}\,s_{{2}} \right)s_{n}({\bold w})+\left(s_{{1}}+\boldsymbol{2}\,s_{{2}} \right)s_{n-1,1}({\bold w}) +s_{{2}}\,s_{{n-2,2}}({\bold w}) & \forall\, n\geq4.
\end{cases}
\end{equation*}
\end{conj}
When the total degree 3, we have the following conjecture verified up to $n=6$ and $\ell=3$:
\begin{conj}\label{ConjeturaGrado3Familia}
Le $\mathcal{T}_3$ be the family of all monomials of total degree 3 in the variables $x_{11},\ldots,x_{1n}$. Then the graded Frobenius characteristic of the associated polarization module is given by table
 \ref{table:FrobeniusGrado3Familia}
\end{conj}
\begin{table}[H]
\caption{Frobenius series for the family of monomials of degree 3}
\label{table:FrobeniusGrado3Familia}
\centering
\scalebox{1.0}{
\begin{tabular}{|c|c|}
\hline
\multirow{3}{*}{$(1+s_1+s_2+s_3)s_{1}({\bold w})$} & \\&$n=1$\\& \\
\hline
\multirow{3}{*}{$(1+s_1+\boldsymbol{2}\,s_2+\boldsymbol{2}\,s_3)s_{2}({\bold w})
+(s_1+s_2+s_{1,1}+\boldsymbol{2}\,s_{3})s_{1,1}({\bold w})$} & \\&$n=2$\\& \\
\hline &\\
$(1+s_1+\boldsymbol{2}\,s_2+\boldsymbol{3}\,s_3)s_{3}({\bold w})
+(s_1+\boldsymbol{2}\,s_2+s_{1,1}+\boldsymbol{3}\,s_3)s_{2,1}({\bold w})$ &\\
$+(s_{1,1}+s_3)s_{1,1,1}({\bold w})$ & \\ & \multirow{-4}{*}{$n=3$} \\
\hline
&\\ 
$\left(1+s_{{1}}+\boldsymbol{2}\,s_{{2}}+\boldsymbol{3}\,s_3\right)
s_{4}({\bold w})+\left(s_{{1}}+\boldsymbol{2}\,s_{{2}}+s_{1,1}
+\boldsymbol{4}s_{3}\right)s_{3,1}({\bold w})$ & \\
$+(s_{{2}}+s_3)s_{{2,2}}({\bold w})+(s_{1,1}+s_3)s_{2,1,1}({\bold w})$ & \\ &
\multirow{-4}{*}{$n=4$} \\
\hline
&\\
$\left( 1+s_{{1}}+\boldsymbol{2}\,s_{{2}}+\boldsymbol{3}\,s_3\right)
s_{5}({\bold w})+\left(s_{{1}}+\boldsymbol{2}\,s_{{2}}+s_{1,1}
+\boldsymbol{4}s_{3}\right)s_{4,1}({\bold w})$ &\\
$+(s_{{2}}+\boldsymbol{2}\,s_3)s_{{3,2}}({\bold w})
+(s_{1,1}+s_3)s_{3,1,1}({\bold w})$ &\\ & \multirow{-4}{*}{$n=5$} \\
\hline
&\\
$\left(1+s_{{1}}+\boldsymbol{2}\,s_{{2}}+\boldsymbol{3}\,s_3\right)
s_{n}({\bold w})+\left(s_{{1}}+\boldsymbol{2}\,s_{{2}}+s_{1,1}+\boldsymbol{4}s_{3}\right)s_{n-1,1}({\bold w})$ &\\
$+(s_{{2}}+\boldsymbol{2}\,s_3)s_{{n-2,2}}({\bold w})
+(s_{1,1}+s_3)s_{n-2,1,1}({\bold w})
+s_{3}\,s_{n-3,3}({\bold w})$&\\ & \multirow{-3}{*}{$\forall\,n\geq 6.$}\\
\hline
\end{tabular}} 
\end{table}
\noindent Corollary \ref{CorolarioGrado2}, Corollary \ref{CorolarioGrado3} and independent verifications with Maple and Sage, lead us to conjecture the following affirmations about the Frobenius characteristic of $\mathcal{M}_f$ when $f$ is homogeneous of degree 4 or 5:
\begin{conj}
The classification given by tables  \ref{table:FrobeniusGrado4} and \ref{table:FrobeniusGrado5} is complete (up to $n$-exceptions), that is, if $f$ is any homogeneous diagonally symmetric polynomial of degree 4 (respectively, degree 5) then the Frobenius characteristic of the module $\mathcal{M}_f$ is one of the formulas in the table \ref{table:FrobeniusGrado4} (respectively, table\ref{table:FrobeniusGrado5}).
\end{conj}
Looking at tables \ref{table:FrobeniusGrado4} and \ref{table:FrobeniusGrado5} we are lead to think that
\begin{conj}
If $f$ is an homogeneous symmetric polynomial of degree ${\bold d}$. Then, there exist a monomorphism $\varphi:\mathcal{M}_{f}\longrightarrow\mathcal{M}_{h_{m}}$.
\end{conj}
A description of Hilbert series of polarization modules and further research directions are available in \cite{Blandin}. The next conjecture is already settled for degree 2 and 3, in higher degree, we believe that 
\begin{conj}\label{Hpositividad}
Let $f$ be any homogeneous symmetric polynomial. The Hilbert series of the module $\mathcal{M}_f$ is $h$-positive, that is, there are $a_{\mu}\in\naturals$ such that $\mathcal{M}_{f}({\bold q})=\sum_{\mu}a_{\mu}h_{\mu}({\bold q}).$
where the sum runs over the set of partitions $\mu$ of integers less or equal to $\deg(f)$.
\end{conj}
\section{Partial results}
\begin{table}[H]
\caption{Frobenius characteristic for degree 4}
\label{table:FrobeniusGrado4}
\centering
\scalebox{.8999}{
\begin{tabular}{|c|c|}
\hline
\multirow{3}{*}{$(1+s_1+s_2+s_3+s_4)s_{n}({\bold w})$} &
 \\ & $p_{1}^4$\\ & \\
\hline
\multirow{3}{*}{$(1+s_1+s_2+s_3+s_4)s_{n}({\bold w})
+(s_1+s_2+s_3)s_{{n-1,1}}({\bold w})$} & \\ & $p_{4}$ \\ & \\
\hline
\multirow{3}{*}{$(1+s_1+s_2+s_3+s_4)s_{n}({\bold w})
+(s_1+s_2+s_3)s_{{n-1,1}}({\bold w})+s_2 s_{{n-2,2}}({\bold w})$} &
\\& $e_4$\\ & \\
\hline
\multirow{3}{*}{$(1+s_1+\boldsymbol{2}s_2+\boldsymbol{2}s_3+s_4)s_{n}({\bold w})
+(s_1+\boldsymbol{2}s_2+s_3)s_{{n-1,1}}({\bold w})$} & \\ &$e_{_{31}}$\\ & \\
\hline
\multirow{3}{*}{$(1+s_1+\boldsymbol{2} s_2+\boldsymbol{2} s_3+s_4) s_{n}({\bold w})
+(s_1+\boldsymbol{2} s_2+s_3) s_{{n-1,1}}({\bold w})
+s_2 s_{{n-2,2}}({\bold w})$} &
$s_{_{211}}$\\ & $h_{_{22}}$\\ & $m_{_{211}}$ \\
\hline
\multirow{3}{*}{$(1+s_1+\boldsymbol{2} s_2+\boldsymbol{2} s_3+s_{21}
+s_4)s_{n}({\bold w})
+(s_1+s_2+s_{11}+s_3) s_{{n-1,1}}({\bold w})$} &
$p_{_{211}}$ \\ & $e_{_{211}}$ \\ & $h_{_{211}}$\\
\hline
\multirow{3}{*}{$(1+s_1+\boldsymbol{2} s_2+\boldsymbol{2} s_3
+s_{21}+s_4) s_{{n}}({\bold w})
+(s_1+\boldsymbol{2} s_2+s_{11}+s_3) s_{{n-1,1}}({\bold w})$} &
$h_{_{31}}$ \\ & $m_{_{31}}$ \\ & $p_{_{31}}$\\
\hline
\multirow{3}{*}{$(1+s_1+\boldsymbol{2} s_2+s_3+s_{21}+s_4) s_{n}({\bold w})
+(s_1+\boldsymbol{2} s_2+s_{11}+s_3)s_{{n-1,1}}({\bold w})
+s_{2} s_{{n-2,2}}({\bold w})$} &
\\ & $m_{_{22}}$ \\ & \\
\hline
\multirow{5}{*}{$(1+s_1+\boldsymbol{2} s_2+\boldsymbol{2} s_3+s_{21}+s_4) s_{n}({\bold w})
+(s_1+\boldsymbol{2} s_2+s_{11}+s_3) s_{{n-1,1}}({\bold w})+s_2 s_{{n-2,2}}({\bold w})$} &
$s_{_{4}}$ \\ & $s_{_{31}}$ \\ & $s_{_{22}}$ \\ & $e_{_{22}}$ \\ & $p_{_{22}}$\\
\hline
\end{tabular}}
\end{table}

\begin{table}[ht!]
\caption{Frobenius characteristic for degree 5}
\label{table:FrobeniusGrado5}
{\footnotesize
\begin{tabular}{|c|l|}
\hline
\multirow{3}{*}{$(1+s_1+s_2+s_3+s_4+s_5)s_{n}({\bold w})$}
& \\
&$p_{1}^5$\\
&\\
\hline
\multirow{3}{*}{$(1+s_1+s_2+s_3+s_4+s_5)s_n({\bold w})
+(s_1+s_2+s_3+s_4)s_{n-1,1}({\bold w})$} & \\
& $p_{5}$ \\
& \\
\hline
\multirow{3}{*}{$(1+s_1+s_2+s_3+s_4+s_5)s_n({\bold w})
+(s_1+s_2+s_3+s_4)s_{n-1,1}({\bold w})+(s_2+s_3)s_{n-2,2}({\bold w})$} &  \\ & $e_5$ \\ &  \\
\hline
\multirow{3}{*}{
$(1+s_1+\boldsymbol{2}s_2+\boldsymbol{2}s_3+\boldsymbol{2}s_4+s_5)s_n({\bold w})
+(s_1+\boldsymbol{2}s_2+\boldsymbol{2}s_3+s_4)s_{n-1,1}({\bold w})
+(s_2+s_3)s_{n-2,2}({\bold w})$}
& $m_{_{2111}}$ \\
& $s_{_{2111}}$\\
&$e_{41}$\\
\hline
\multirow{3}{*}{$\begin{array}{c}
(1+s_1+\boldsymbol2s_2+\boldsymbol2s_3+s_{21}+\boldsymbol2s_4+s_{31}+s_5)s_n({\bold w})
+(s_1+\boldsymbol2s_2+s_{11}+\boldsymbol2s_3+s_{21}+s_4)s_{n-1,1}({\bold w})\\
+(s_2+s_3)s_{n-2,2}({\bold w})
\end{array}$}
& \\
& $s_{221}$ \\
& \\
\hline
\multirow{2}{*}{$(1+s_1+\boldsymbol{2}s_2+\boldsymbol{2}s_3+s_{21}+\boldsymbol{2}s_4+s_{31}+s_5)s_n({\bold w})
+(s_1+\boldsymbol{2}s_2+s_{11}+\boldsymbol{2}s_3+s_{21}+s_4)s_{n-1,1}({\bold w})$}
& $m_{41}$ \\
& $p_{41}$\\
\hline
\multirow{9}{*}{$\begin{array}{c}
(1+s_1+\boldsymbol{2}s_2+\boldsymbol{3}s_3+s_{21}+\boldsymbol{2}s_4+s_{31}+s_5)s_n({\bold w})
+(s_1+\boldsymbol{2}s_2+s_{11}+\boldsymbol{3}s_3+s_{21}+s_4)s_{n-1,1}({\bold w}) \\
+(s_2+s_3)s_{n-2,2}({\bold w})\end{array}$}
&$h_5$ \\ & $h_{41}$\\
&$h_{32}$ \\ & $h_{221}$\\
&$p_{221}$ \\ & $s_{41}$\\
&$s_{32}$ \\ & $s_{311}$\\
&$e_{221}$ \\ & $m_{311}$\\
\hline
\multirow{4}{*}{$\begin{array}{c}
(1+s_1+\boldsymbol{2}s_2+\boldsymbol{2}s_3+s_{21}+\boldsymbol{2}s_4+s_{31}+s_5)s_n({\bold w})
+(s_1+\boldsymbol{2}s_2+s_{11}+\boldsymbol{3}s_3+s_{21}+s_4)s_{n-1,1}({\bold w})\\
+(s_2+s_3)s_{n-2,2}({\bold w})
\end{array}$}
&$p_{32}$ \\ & $e_{32}$ \\
&$m_{32}$ \\ & $m_{221}$\\
\hline
\multirow{3}{*}{
$(1+s_1+\boldsymbol{2}s_2+\boldsymbol{2}s_3+s_{21}+\boldsymbol{2}s_4+s_{31}+s_5)s_n({\bold w})
+(s_1+s_2+s_{11}+s_3+s_{21}+s_4)s_{n-1,1}({\bold w})$}
& $p_{2111}$\\
& $h_{2111}$\\
& $e_{2111}$\\
\hline
\multirow{3}{*}{
$(1+s_1+\boldsymbol{2}s_2+\boldsymbol{3}s_3+s_{21}+\boldsymbol{2}s_4+s_{31}+s_5)s_n({\bold w})
+(s_1+\boldsymbol{2}s_2+s_{11}+\boldsymbol{2}s_3+s_{21}+s_4)s_{n-1,1}({\bold w})$}
&$e_{311}$ \\
&$h_{311}$ \\
&$p_{311}$ \\
\hline
\end{tabular}}
\end{table}

\clearpage

\acknowledgements\label{Agradecimientos}

Thanks to Fran\c cois Bergeron and Franco Saliola for their advice  during my Ph.D. studies. Thanks to Franco Saliola and Eduardo Blazek for helpful discussions and many independent computations Sage. Thanks to Adolfo Rodr\'iguez for his help with Maple. Thanks to Yannic Vargas, Alejandro Morales, Marco P\'erez and Luis Serrano for many valuable suggestions to reach the final version of this paper. To my fiancee Oana-Andreea Kosztan for love and support with \LaTeX\ during my Ph.D. Thesis. I want to thank Francois Bergeron for showing me this beautiful subject.

\nocite{*}

\bibliographystyle{abbrvnat}


\bibliography{bibliografia}

\end{document}